\documentclass[12pt,letterpaper]{amsart}
\usepackage{verbatim}
\usepackage{amsmath}
\usepackage{amssymb}
\begin{document}
\author{David W. Farmer}
\thanks{The REUs described in this paper were supported by
AIM and by grants from the NSF, including the Focused Research
Group grant DMS 0244660.}
\address{American Institute of Mathematics}
\email{farmer@aimath.org}
\title{The AIM REU: individual projects with a common theme}
\maketitle

The principle behind the
American Institute of Mathematics (AIM) REU is that mathematicians with an
active research program have ideas for more projects than
they have time to work on.  
By enlisting the help of enthusiastic undergraduates, it is possible
to broaden the range of problems in which they
are actively involved.

The AIM REU is run by David Farmer and Brian Conrey, along with
one or two of their current postdocs.  Conrey and Farmer, and 
their postdocs, have research interests
in analytic number theory: the Riemann zeta-function, modular forms,
random matrix theory, elliptic curves, etc.  The REU projects all
come directly from their ongoing research programs, so there 
are natural relationships between the projects.
The fact that there is an underlying theme between the
various student projects is one
of the strengths of the program.  
There is a wide variety of interaction among the students and the advisers,
with all the students developing an interest in each other's projects,
and all the advisers being able to work with all the students.
The students become part of an active research community.

\section{Key features of the AIM REU}

\textbf{Individual projects. } Every student chooses their own
individual research problem.  
This allows the student to take full responsibility
for their own project and 
more closely models traditional mathematics research.

The students all work on related projects, although it often
takes them many weeks to understand the connections!  Through
weekly progress reports and a variety of daily interactions
with the advisers and other students, the participants spend
a lot of time talking to each other, explaining their work,
and learning about each other both on a mathematical and a personal
level.  

By having all the projects in related areas we are able
to foster a sense of community where everyone feels a part of a larger
endeavor.  This allows us to have all the students work on separate
projects without the likelihood of feeling isolated.

\textbf{Immediate work.  }  Students choose their problem on the first
day and begin working immediately, without any preliminary background
reading.   This feature is somewhat unusual, considering that the
problems are in an area that is considered to require a reasonable
amount of background.  We elaborate on this in the next two sections of
this paper.

\textbf{Ongoing research. }
The student projects arise naturally from the research programs
of the mathematicians running the REU.
One consequence is that the specific
projects cannot be finalized until the last minute.  While it is possible
to give a general idea of the program (students usually want to
know the possible projects and how their particular project
will be decided), we always have to consider recent
developments before the beginning of the REU.

\textbf{Postdoctoral training. }  
The postdocs are involved in all aspects of the REU: everything
from selecting appropriate projects to helping the students
prepare their talks in the final week.
This is valuable training because there is an increasing
expectation for junior faculty to advise undergraduate projects
(especially at liberal arts colleges). 

This ``training'' begins shortly before the program when the list of possible
projects is developed.  It is not required that the postdocs contribute
a problem (although they usually do).  The senior mathematicians 
and postdocs
have 2 or 3 meetings which start with a discussion of ideas
for projects.  The senior mathematicians examine all aspects of
the proposed problems, being as explicit as possible about 
which properties make a problem suitable, or not, for the REU.
By the final meeting, the discussions turns to the critical
issue of how to introduce the problems to the student so that they
can make a reasonable choice for their project and begin working
immediately.  This process is not easy to capture on paper,
so it is a big help for the postdocs to see it in action.

Regular meetings between the directors and the postdocs occur
during the program.  These are used to monitor student progress
and also to discuss issues
such as:  what is the right amount of help to give a student
who is ``stuck,'' how to prepare a student to give a talk,
etc.  These meetings are in addition to the daily interaction of
the mentors and the students, which gives the postdocs an opportunity
to observe the more experienced mentors and to call on help as needed.

\textbf{Mathematica. } Many of our REU projects 
begin with experimentation in \textsl{Mathematica},
and in many cases the work is done primarily with the
help of \textsl{Mathematica} and other computer-algebra
systems.

\textbf{Standard topics. }  The AIM REU has many facets which
feature in most REUs, but we mention them for completeness.
Throughout the program there are lectures on material related to
the student projects.  All the students learn LaTeX and write
a final report in the style of a research paper.  Also, the students
give a formal talk on their work and are encouraged to give a talk
when they return to their home institution.  The mentors and the
students all see each other several times every day and often eat
lunch together.  Each student has their own desk in a shared office.

\section{Identifying appropriate projects}

Prior to the REU we must choose projects which allow our
students to select a problem on the first day and begin
work immediately, after only a few verbal preliminaries
and without any background reading.  We have identified three
key criteria for deciding if a problem is appropriate for this approach,
and one ``non-criterion'' which is best ignored.

\textbf{1. How would I start to work on the problem? }
The problem is part of my research, so obviously I find it interesting.
Suppose I had time to work on the problem: 
\emph{what would be my first step? }
If I don't have an answer to that
question, then clearly the problem is not appropriate
for a student.   On the other hand, if I have a good
idea of what I would do, then it is possible that the
problem may be good for the REU.

\textbf{2. Could an undergraduate execute my first steps on the problem? }
I know what I would try if I were to work on the problem.  Would
an undergraduate be able to do that work?  For example,
if I would start with some
numerical experiments in \textsl{Mathematica}, or I would start by considering
the problem for some low-degree polynomials, or if I would enumerate
the first few cases by brute-force, then perhaps a student would
be able to do it, too.

Note: this question only asks if the student could do the work,
at this point we don't consider whether or not the student could
understand \emph{why} that work relates to the problem.

\textbf{3. It there something related to the problem which the student
could begin doing immediately? }
An effective way to engage a student on a problem is to have them
begin work on the first day: work which will lead them to see
something interesting.
This only needs to be in the same general area as the
intended problem, because its purpose is to quickly put the
student into a research-like environment and to start them towards
an understanding of the problem.  Once they have seen something
new (to them), they want to understand their
observations (which motivates them to learn background material),
and they will want to keep going (which motivates them to get to the
real problem).

For example, I have had several students do research on modular forms.
A good first task for them is:  use \textsl{Mathematica} to expand the
following infinite product into a sum:
\begin{equation}
q \prod_{j=1}^\infty (1-q^j)^{24} = \sum_{n=1}^\infty \tau(n) q^n .
\end{equation}
It doesn't take much knowledge of \textsl{Mathematica} to truncate the product
on the left and expand it into a sum, which gives the first several
values of $\tau(n)$.  The student will ``discover'' that
$\tau(2)\tau(3)=\tau(6)$ and $\tau(2)\tau(5)=\tau(10)$, and in
general $\tau(n)\tau(m)=\tau(nm)$ if $n$ and $m$ are relatively
prime.  The proof of this result is accessible to an undergraduate
(although it is not found in the undergraduate curriculum),
and the student is motivated to see the proof.

\textbf{4. You don't need to fill in the gap. }
After the student has started, you still have to chart a course
to their chosen problem.  While it is helpful to have an idea
of this before the REU begins, I have not found this to be critical.
You will work together with the student to find a path from the
initial task to the problem.  You have already determined that both
ends of that path are accessible to the student, so you have to
trust that you can connect them in an understandable way.
Often you can identify a huge chasm which seems hopeless for
the student to breach.  In these cases a few carefully chosen 
``black boxes'' (which the student will have to accept on faith)
can be helpful.  Your goal is not to give an entire graduate
course in the subject, but to put them in a position to work on their
problem.

It takes a lot of work to bring a student up to speed on a problem,
and this is best done in partnership with them because you need
to take into account their particular strengths and weaknesses.

\section{The beginning of the project}

Ten weeks (or less) is a very short amount of time 
to accomplish meaningful research, so it is important for students
to begin right away.   In the AIM program a variety of possible
projects are described on the first morning.  Shortly after
lunch each student chooses a problem and begins work.
Students are never asked to read or learn anything before
beginning their project, and the students are not formally introduced
to any of the 
background of their project until their work is well underway.

\textbf{Real work begins the first day. }  This is the aspect of
the program which many people find surprising, including the
students.  The first morning all the mentors describe
possible projects.  Technical terms are suppressed as much as
possible, and the emphasis is on the \emph{methods} that the
project will involve.   For example, some projects involve a lot of
computer work, others have a more geometric flavor, etc.
By lunchtime all the projects have been described, and the students
are told that after lunch they can ask as many questions as they want,
after which they must choose a project and begin working
immediately.  They are not forced to stay with that project
if they don't like it, but they have to give it a serious try.
Should two people want to do the same project, 
the mentors can usually find two independent sub-projects
with the same theme.
The mentors and students eat lunch together and then return to AIM.

The after-lunch meeting usually takes an hour, and after the student 
questions are answered we leave them alone for some time to talk
it over.   Some students are uncomfortable with this process, and many
explain that they would like to learn more about the problems
before choosing.  We tell them: ``You have spent years learning 
mathematics, and you will have plenty of time to
learn more mathematics during the semester.  
If you spend too much time learning the background, then
you won't have time to  \emph{do} the  mathematics.''

\textbf{Starting a project. }  
How are they are supposed to
start working without any background?  
The answer is to give them specific tasks which they can 
do using only what they already know, and which bear some relation
to their project.  The tasks should lead them to see some
interesting phenomenon.  This could involve a pen-and-paper
calculation, or a computer calculation (we use \textsl{Mathematica} a lot),
or drawing a picture, etc.
There is no need to explain the connection between the task and
the project -- that will come later.

The students complete their first task, which 
hopefully leads them to an interesting observation.  Now they
are hooked.  They want to move on to what is next, they want
to know the relationship to their project, they want to
understand the background, etc.  Of course the student will
have to learn some background at some point, but that can
be done in parallel with the ``real'' work, and the research provides
the motivation.  You have asked them to trust you and begin
without knowing where they are going, but once they get started
you have to provide them the means to fill in the gaps.

I intend this approach to be an explicit rejection of the following
method, which has no place in an REU:
``Go read
these books and papers and when you understand them come back and I'll
tell you about your project.''

I do not claim that every research project can be started in this way.
But a very large number of those \emph{which are suitable for
REU projects} can indeed fit our model.  It takes some time and creativity,
and this is a key part of the process of identifying suitable projects.

\textbf{Students have individual projects. }  Each student has
his or her own separate project.  This project is under the primary
supervision of one of the mathematicians, although at some point
during the summer every mathematician works with every student. 
This helps expose the student to different perspectives and
also makes it possible for the mathematicians to go to conferences
without compromising the attention given to the students.

\section{Miscellaneous}

Some features of the AIM REU which we recommend for other programs:

\textbf{Conference call prior to the summer.  }
I arrange a conference call with all the students who have an offer
to join our REU.  This is efficient because there are many
questions which nearly
all students have, so I can answer them all at once.  It also
helps the students to find out a little about each other and to begin
building a sense of community.

\textbf{Weekly reports. }
Every Monday morning all the students give progress reports.
We specify that the reports should be self-contained and should not assume
the audience recalls the details of their previous report.
What is their project and why is it interesting?  
What did they do last week?
What is the overall plan and what specifically will they do
this week?

We meet with the students before the report to discuss what they
will say.  Usually they need encouragement to be expansive about the
``big picture.''  We meet with them again after the report to offer
constructive criticism, usually phrased in the form: 
``Next week you might try...''

Initially the students do not like giving these reports, but in
our exit surveys they tell us that the reports were valuable and
they are glad they did them.
In addition to making them more comfortable about giving a formal
presentation of mathematics,
they specifically indicate that the reports help them to keep
their project in perspective and to plan their week.

\textbf{LaTeX at the first opportunity. }  Have the students
learn LaTeX the first time they do any work that could be
part of their final report.  This could be original work,
or it could be some standard calculation
that they worked through which could be part
of the introduction in their paper.  Give them a simple
template file, and just have them type something small.  Don't
ask for a formal writeup or introduction:  just have them
start learning how to LaTeX.  Encourage them to keep writing
up anything that could be useful.

\textbf{Preparing final talks. }  We host a mini-symposium during
the final two days of the program
where all the students give 40-50 minute
talks.  We treat this just like a session at a conference.

Sally Koutsoliotas and I have developed an effective approach
to preparing students to give talks.  We prepared
a guide, available at
\texttt{aimath.org/mathcommunity/}. 
Here is a summary:
it is based on a series of four meetings.

The first meeting is just
a conversation with the student.  We identify the main point
of their talk and work backward to determine the absolute minimum
amount of information needed to \emph{understand} and
\emph{appreciate} their main point.  The student goes away and
prepares a rough outline of their talk, usually in the form of 
handwritten drafts of slides on paper. (It is critical that the
student prepare something very rough, because it will probably have to
be completely taken apart and you don't want them to resist your
suggested changes).

The second meeting begins with the student summarizing their
main point, \emph{and then} presenting their rough draft by
laying it all out on a table.  The mentor and student work
together to identify the
different parts of the talk, stressing the importance of the
\emph{transition} between each part, and the \emph{flow}
of ideas.  The student goes away to make a complete version
of their slides.

The third meeting involves laying out the slides on a big table
(use printouts if the talk is on the computer) and grouping them
according to the different parts of the talk.  We evaluate
the total number of slides, the content and layout of each slide,
and the all-important transition between the different parts of the talk.
The student goes away and makes final versions of the slides.

The fourth meeting is the first time the student actually gives
a practice talk.  
Tell them to just talk and don't worry about how
long it will take (you don't want them rushing at the end: if material
has to be cut, it probably isn't at the end of the talk).
Then sit down together at a table with all the slides laid out,
and try to discuss only the 4 or 5 most important things they need
to change.  If appropriate, that is, if they have to make extensive
revisions,  schedule another practice talk.  If there are minor changes,
such as the talk was slightly too long, you may want to suggest that they
practice in front of a friend.

\textbf{Sample papers. }  To help the students in their final writeup
we prepare a folder of sample papers.  These are research papers from the 
general subject area of the REU and which have good overall organization,
a well-written abstract, and an accessible introduction.   These papers
are models to show the student what their final paper should look like.

\section{The benefits, mathematical and otherwise}

Running an REU is a rewarding experience, particularly when 
you can quickly bring students up to speed and make them part
of your research group.  However, it is a huge amount of work,
and to do it right you have to be willing to dedicate the majority
of your time to it.

Is the total mathematical output
larger than than it would have been if, instead of the REU, I just spent all day
working on my own research?  I think the answer is `yes,' but
in some sense the question is meaningless because I would not have
the stamina to focus all-day-every-day on my research.
The REU definitely increases the breadth of projects I am involved in,
and that makes me a better mathematician.

\end{document}